\documentclass[10pt,draft,oneside]{amsart}

\usepackage{enumerate}
\usepackage{amsthm}
\usepackage{amsfonts}
\usepackage{amsmath}

\newtheorem{theorem}{Theorem}
\newtheorem{proposition}[theorem]{Proposition}
\newtheorem{corollary}[theorem]{Corollary}
\newtheorem{lemma}[theorem]{Lemma}
\theoremstyle{definition}
\newtheorem{remark}[theorem]{Remark}
\newtheorem{question}[theorem]{Question}

\def\Diff{\Delta}
\def\NN{\mathbb{N}}
\def\RR{\mathbb{R}}
\def\ekv{\sim}
\def\eqegy{=}

\begin{document}

\title{Decomposition as the sum
of invariant functions with respect to commuting transformations}

\author{B\'{a}lint Farkas and Szil\'{a}rd Gy. R\'{e}v\'{e}sz}
\date{\today}
\maketitle

\begin{abstract} As natural generalization of various investigations
in different function spaces, we study the following problem. Let
$A$ be an arbitrary non-empty set, and $T_j$ ($j=1,\dots,n$) be
arbitrary commuting mappings from $A$ into $A$. Under what
conditions can we state that a function $f:A\to \RR$ is the sum of
``periodic'', that is, $T_j$-invariant functions $f_j$? (A function
$g$ is periodic or invariant mod $T_j$, if $g\circ T_j=g$.) An
obvious necessary condition is that the corresponding multiple
difference operator annihilates $f$, i.e.,
$\Delta_{T_1}\dots\Delta_{T_n}f= 0$, where $\Delta_{T_j}f:=f\circ
T_j-f$. However, in general this condition is not sufficient, and
our goal is to complement this basic condition with others, so that
the set of conditions will be both necessary and sufficient.

\vskip1em \noindent{\small \textbf{Mathematics Subject
Classification (2000):} Primary 39A10. Secondary 39B52,
39B72.\\[1em]
} \vskip1em
\end{abstract}
\section{Introduction}
Let $f:\RR\to \RR$ be a function which is a sum of finitely many periodic functions
\begin{equation}\label{eq:periodd}
f=f_1+f_2+\cdots +f_n,\qquad f_i(x+\alpha_i)=f_i(x)\quad\mbox{for all $x\in\RR$, $i=1,\dots,n$}
\end{equation}
with some fixed numbers $\alpha_i\in\RR$. For $\alpha\in\RR$ let $\Diff_\alpha$ denote the
(forward) difference operator
\begin{equation*}
\Diff_\alpha:\RR^\RR\to\RR^\RR,\quad\Diff_\alpha g(x):=g(x+\alpha)-g(x)~.
\end{equation*}
Then the $\alpha_i$-periodicity of $f_i$ above means $\Diff_{\alpha_i}f_i\eqegy0$, and because
the difference operators are commuting, we also have that
\begin{equation}\label{eq:kern}
\Diff_{\alpha_1}\Diff_{\alpha_2}\dots\Diff_{\alpha_n}f\eqegy 0~.
\end{equation}
The starting point of this work is the problem, if the converse of
the above statement is also true, i.e., does \eqref{eq:kern} imply
the existence of a periodic decomposition \eqref{eq:periodd}?

Naturally, this question can be posed in any given function class
$X\subset \RR^{\RR}$, when we have $f\in X$ and \eqref{eq:kern} and
want a decomposition \eqref{eq:periodd} within the class, i.e., we
require also $f_i\in X$ ($i=1,\dots,n$). If the answer is
affirmative, the class $X$ is said to have the \emph{decomposition
property}. It is easy to see that $\RR^\RR$ or $C(\RR)$ (space of
continuous functions) do \emph{not} have this property. Indeed, let
$n=2$ and $\alpha_1=\alpha_2=\alpha$. The identity function
$f(x):=x$ satisfies $\Delta_{\alpha} \Delta_{\alpha} f \eqegy 0$,
but $f(x)=x$ fails to be the sum of two $\alpha$-periodic functions
(as then $f$ would be periodic itself). This shows that the
implication \eqref{eq:kern} $\Rightarrow$ \eqref{eq:periodd} fails
even in the simplest possible case, and that further conditions
either on the transformations or on the functions are needed in
general.

The study of such decomposition problems originated from
I.~Z.~Ruzsa. In particular, he showed that the identity function
$f:\RR\to\RR$, $f(x)=x$ has a decomposition into the sum of
$\alpha$- respectively $\beta$-periodic functions provided that
$\alpha$, $\beta$ are incommensurable. Later M.~Wierdl observed that
in the space of arbitrary real-valued functions the difference
equation \eqref{eq:kern} implies \eqref{eq:periodd} if the steps
$\alpha_i$ are linearly independent over $\mathbb{Q}$ (see
\cite[Lemma, p.~109]{wierdl:1984}). We extend this result in
Corollary \ref{cr:irrational} to the case when the shifts $\alpha_i$
are only assumed to be pairwise incommensurable.

 Answering also a question of M.~Laczkovich, Wierdl \cite{wierdl:1984} showed that for $n=2$ identity \eqref{eq:kern} implies \eqref{eq:periodd} in the space
 $BC(\RR)$ (bounded continuous functions). The proof for general $n$, i.e., that $BC(\RR)$ has the decomposition property was later done by M.~Laczkovich and Sz.~Gy.~R\'{e}v\'{e}sz
\cite{laczkovich/revesz:1989}. Moreover, Laczkovich and R\'{e}v\'{e}sz
generalized the problem to many other function classes in fact
considering the derived problem in topological vector spaces
\cite{laczkovich/revesz:1990}. Later Z.~Gajda \cite{gajda:1992} gave
 alternative proofs based on Banach limits that the spaces $B(\RR)$ (bounded functions) and $UCB(\RR)$ (uniformly
continuous and bounded functions) have the decomposition property.
Recently, interesting results and examples were found by
V.~M.~Kadets and S.~B.~Shumyatskiy \cite{kadets/shumyatskiy:2000,
kadets/shumyatskiy:2001} about the decomposition problem in various
Banach spaces.

In a different direction, T.~Keleti \cite{keleti:1996, keleti:1997,
keleti:1998} studied related problems and was led to a negative
answer regarding the existence of a continuous periodic
decomposition of continuous functions on $\RR$ provided only that a
measurable decomposition exits on $\RR$, see \cite[Theorem
4.8]{keleti:1996}.

In the present paper we do not restrict to any particular function class, and we neither assume
any particular structural properties like smoothness etc.~of the transformations. The present
work attacks the decomposition problem in the whole space of functions $\RR^A$ with respect to
arbitrary commuting operators in $A^A$.

\vskip1em

Throughout this note $A$ will denote a fixed nonempty set. We will
consider various self maps $T:A\to A$, called
\emph{transformations}, and to such a transformation we associate
the corresponding \emph{shift operator} $T$ (denoted by the same
symbol) as $T(f):=f\circ T$ and the $T$-\emph{difference operator}
$\Diff_T:\RR^A\to \RR^A$ defined as
$$
(\Diff_T f):=  T(f) - f~,\qquad\qquad (\Diff_T f)(x)=f(Tx)-f(x)~.
$$
A function $f$ satisfying $\Diff_T f\eqegy 0$ is called \emph{$T$-invariant}.

A \emph{$(T_1,\dots,T_n)$-invariant decomposition} of some function $f$ is a representation
\begin{equation}\label{invaridecdef}
f=f_1+f_2+\cdots +f_n ~,\qquad\text{where}\qquad  \Diff_{T_j} f_j \eqegy0 \quad (j=1,\dots,n)~.
\end{equation}
As mentioned above, we do not assume any properties like smoothness,
boundedness, injectivity or surjectivity etc., neither on the
transformations nor on the functions, except that all the occurring
transformations and functions are defined over the whole set $A$ and
that the occurring transformations must commute. For pairwise
commuting transformations $T_i$ the functional equation
\begin{equation}\label{eqdifference}
\Diff_{T_1}\dots \Diff_{T_n} f = 0
\end{equation}
evidently implies for every $k_1,\dots,k_n\in\NN$ (where also $0\in\NN$ in our present
terminology) the equation
\begin{equation}\label{eqdifferencekj}
\Diff_{T_1^{k_1}}\dots \Diff_{T_n^{k_n}} f = 0~.
\end{equation}
Now in this general setting our basic question sounds: Does the functional equation
\eqref{eqdifference} imply the existence of some $(T_1,\dots,T_n)$-invariant decomposition
\eqref{invaridecdef}? As mentioned above, in general this is not the case. Therefore, we look
for further conditions, which, together with \eqref{eqdifference}, are not only necessary, but
also sufficient to ensure that such an invariant decomposition exists. More precisely, in the
next section we focus on complementary conditions -- functional equations -- on the functions,
which they must satisfy in case of existence of an invariant decomposition \eqref{invaridecdef}
and which equations will also imply existence of such a decomposition. In the third section we
define a further, still quite general property of transformations, implying that the difference
equation \eqref{eqdifference} also suffices for the existence of an invariant decomposition
\eqref{invaridecdef}.

With this general framework the pure combinatorial nature of the
problem is quite apparent. Since similar questions arise quite often
in various settings, the present formulation may help understanding
some related problems as well. To emphasize the combinatorial
structure, one may reformulate the whole problem so as to consider
$A$ as the vertices of a directed and colored graph, with $T_j$
being the set of directed edges, colored by the $j^{\text{th}}$
color. Transformations are defined uniquely on $A$, which means that
the out-degree of any color is exactly $1$ at all points of $A$. The
pairwise commutativity assumption then means that starting out from
a given point and traveling along one blue and one red edge, we
arrive at the same point independently of the order we chose of
these colors. Looking for (color-) invariant functions is the search
of $f_j$ which assume equal values on points connected by a directed
path of the $j^{\text{th}}$ color. Mentioning this interpretation
may reveal the combinatorial nature of the question, although we do
not emphasize this language any longer.

\section{Results for arbitrary transformations}

\begin{theorem}\label{thm:twodecomp}
Let $A$ be a nonempty set and $S,T:A\to A$ commuting transformations, $f\in\RR^A$. The
following are equivalent
\begin{enumerate}[i)]
\item \label{enum:1:i} There exists a decomposition $f=g+h$, with $g$ and $h$ being
$S$- and $T$-invariant, respectively.
\item \label{enum:1:ii}$\Diff_S\Diff_T f\eqegy 0$, and if for some $x \in A$ and
$k,n,k',n'\in\NN$ the equality
\begin{equation}\label{eq:mixingdef}
T^k S^n x=T^{k'} S^{n'}x
\end{equation}
holds, then
\begin{equation}\label{eqcompatibility}
f(T^k x)=f(T^{k'} x)
\end{equation}
must also be satisfied.
\end{enumerate}
\end{theorem}
\begin{proof}

i) $\Rightarrow$ ii):\hskip1em The first part is obvious. Indeed, as $T$ and $S$ commute,
$$\Diff_S\Diff_T f=\Diff_S\Diff_T g+\Diff_S\Diff_T h=\Diff_T\Diff_S g+\Diff_S\Diff_T
h=0+0=0~.
$$
Suppose now that \eqref{eq:mixingdef} holds for some $x \in A$ and
$k,n,k',n'\in\NN$. Then using commutativity of $S$ and $T$, the
$S$-invariance of $g$ and the $T$-invariance of $h$
\begin{align*}
f(T^k x)& =g(T^k x)+h(T^k x)=g(S^nT^k x)+h(x)\\
& = g(S^{n'}T^{k'} x)+h(T^{k'}x)=g(T^{k'}x)+h(T^{k'}x)=f(T^{k'}x)
\end{align*}
follows.

 \vskip1em\noindent ii) $\Rightarrow$
i):\hskip1em Before we can give the decomposition of $f$ we
partition the set $A$. We say that two elements $x,y\in A$ are
equivalent, if for some $k,n,k',n'\in\NN$ the equality
\begin{equation}\label{equivalence}
T^k S^n x=T^{k'} S^{n'}y
\end{equation}
holds. Needless to say that we indeed defined an equivalence
relation $\sim$, hence the set $A$ splits into equivalence classes
$A/\hskip-5pt\sim$, from which by the axiom of choice we choose a
representation system. Obviously it is enough to define $g$ and $h$
on each of these equivalence classes. Indeed, for $x\in A$ the
elements $x$, $Tx$ and $Sx$ are all equivalent, so the invariance of
the desired functions $g,h$ is decided already in the common
equivalence class. So our task is now reduced to defining the
functions $g$ and $h$ on a fixed, but arbitrary equivalence class
$B$. Let $x\in B$ and $x_0$ be the representative of $B$. By
definition, $x\in B$ means $x\sim x_0$, hence the existence of
$k,n,k',n'\in\NN$ satisfying \eqref{equivalence} with $x_0$ in place
of $y$. Set now
\begin{equation}\label{eq:Gdef}
G(n,k,n',k',x):=f(T^{k'}x_0)-f(T^{k}x)+f(x)~.
\end{equation}
Note that here appearance of $n,k,n',k'$ in the argument refers to a
particular combination of powers of $S$ and $T$ showing $x\sim x_0$
rather than arbitrary free variables. First we show that whenever
$l,m,l',m'\in\NN$ provide an alternative relation
\begin{equation} \label{eq:altrel}
T^l S^m x=T^{l'} S^{m'}x_0~,
\end{equation}
then
\begin{equation}\label{eq:Gdefok}
G(n,k,n',k',x)=G(m,l,m',l',x)~.
\end{equation}
By assumption
\begin{equation}\label{twomixedequiv}
T^{k+l'}S^{n+m'}x=T^{k'+l'}S^{n'+m'}x_0=T^{k'+l}S^{n'+m}x
\end{equation}
holds, so using \ref{enum:1:ii}) we obtain $f(T^{k+l'}x)=f(T^{k'+l}x)$. This, together with the
two sides of \eqref{twomixedequiv}, substituted into the definition \eqref{eq:Gdef} of $G$,
yield
\begin{align}\label{eq:Gseq}
&G(n+m',k+l',n'+m',k'+l',x) \notag\\
&\quad=f(T^{k'+l'}x_0)-f(T^{k+l'}x)+f(x)=f(T^{k'+l'}x_0)-f(T^{k'+l}x)+f(x)\notag\\
&\quad=
G(n'+m,k'+l,n'+m',k'+l',x)~.
\end{align}
Again using \eqref{eq:Gdef}, \eqref{twomixedequiv} and the assumption $\Diff_T\Diff_S f\eqegy0$
(in the form that $\Diff_{T^a}f(z)=\Diff_{T^a}f(S^b z)$ for any $a,b\in\NN$) we obtain
\begin{equation*}
\begin{split}
&G(n+m',k +l',n'+m',k'+l',x)-G(n,k,n',k',x)\\
&\quad=\left(f(T^{k'+l'}x_0)-f(T^{k+l'}x)+f(x)\right)-\left(f(T^{k'}x_0)-f(T^{k}x)+f(x)\right)\\
&\quad=\Diff_{T^{l'}}f(T^{k'}x_0)-\Diff_{T^{l'}}f(T^k x)=
\Diff_{T^{l'}}f(T^{k'}S^{n'}x_0)-\Diff_{T^{l'}}f(T^k S^n x)=0~.
\end{split}
\end{equation*}
This shows $G(n+m',k+l',n'+m',k'+l',x)=G(n,k,n',k',x)$, while the same way also
$G(n'+m,k'+l,n'+m',k'+l',x)=G(m,l,m',l',x)$ follows, so now \eqref{eq:Gseq} implies
\eqref{eq:Gdefok}.

All in all, the function
$$
g(x):=G(n,k,n',k',x)
$$
is well defined on $B$ (whence on the whole of $A$). Now $h$ can not
be else than $h:=f-g$. To complete the proof, we show that $g$ and
$h$ have all the necessary properties. Let $x\in A$ and $x_0$ be the
representative of the class $B$ of $x$: for some $n,n',k,k'\in\NN$
we have $T^{k} S^n x=T^{k'} S^{n'}x_0$, and hence also
$$
T^{k} S^n (Tx) =T^{k'+1} S^{n'}x_0\quad\mbox{and~equivalently}\quad T^{k+1} S^n x =T^{k'+1}
S^{n'}x_0~,
$$
so we can write by the definition of $g$
\begin{equation*}
\begin{split}
\Diff_T g(x) =&\bigl(f(T^{k'+1}x_0)-f(T^{k}(Tx))+f(Tx)\bigr)\\
&-\bigl(f(T^{k'+1}x_0)-f(T^{k+1}x)+f(x)\bigr)=\Diff_T f(x)~.
\end{split}
\end{equation*}
As $h:=f-g$, $\Diff_T h \eqegy 0$ is immediate. Finally, we prove that $\Diff_S g \eqegy 0$.
For $x\in B$ we have by \eqref{equivalence} with $x_0=y$, similarly to the above that
\begin{equation*}
\begin{split}
\Diff_S g(x)=&\bigl(f(T^{k'}x_0)-f(T^{k}Sx)+f(Sx)\bigr)\\
&-\bigl(f(T^{k'}x_0)-f(T^k x)+f(x)\bigr)= - \Diff_{T^k} \Diff_S f (x) =0
\end{split}
\end{equation*}
in view of ii).
% \eqref{eqdifferencekj}.
\end{proof}

\begin{remark}
If $TS\neq ST$ then \ref{enum:1:i}) does not imply $\Diff_S\Diff_Tf\eqegy0$ in general.
\end{remark}

\begin{remark}\label{rem:symmass}
Condition \ref{enum:1:i}) is symmetric with respect to the pairs $g,\: S$ and $h,\: T$. This
gives the further equivalent assertion
\begin{enumerate}[i)]\setcounter{enumi}{2}
\item $\Diff_S\Diff_T f\eqegy 0$, and if for some $x \in A$ and
$k,n,k',n'\in\NN$ \eqref{eq:mixingdef} holds, then
$$
f(S^n x)=f(S^{n'} x)
$$
must be satisfied.
\end{enumerate}
\end{remark}

\begin{theorem}\label{th:necessity} Let $T_1,\dots, T_n$ be commuting transformations of $A$ and let $f$ be a real
function on $A$. In order to have a $(T_1,\dots,T_n)$-invariant
decomposition \eqref{invaridecdef} of $f$, the following Condition
(\textasteriskcentered) is necessary. \vskip1em\noindent\hbox
to\hsize{\vtop{\hsize=0.1\hsize\vfill\noindent{\rm
(\textasteriskcentered)}\vfill}\vtop{\hsize=0.9\hsize
\noindent\begin{minipage}{\hsize}

For every $N\leq n$, disjoint $N$-term partition $B_1\cup
B_2\cup\cdots \cup B_N=\{1,2,\dots,n\}$, distinguished elements
$h_j\in B_j$ $(j=1,\dots,N)$, indices $0<k_j, l_j, l_j'\in \NN$,
$(j=1,\dots, N)$ and $z\in A$ once the conditions
\begin{equation}\label{eq:mcond}
T_{h_j}^{k_j} T_{i}^{l_i}z=T_{i}^{l_i'}z\qquad \mbox{for all $i\in
B_j\setminus \{h_j\}$, for all  $j=1,\dots, N$}
\end{equation}
are satisfied, must also
\begin{equation}\label{eq:mconclusion}
\Diff_{T_{h_1}^{k_1}}\dots \Diff_{T_{h_{N}}^{k_N}}f(z)=0
\end{equation}
hold.
\end{minipage}}}\vskip1em
\end{theorem}

\begin{remark}\label{mzerocondi} In case all the blocks $B_j$ are singletons the condition \eqref{eq:mcond} is empty,
so \eqref{eq:mconclusion} expresses exactly \eqref{eqdifferencekj}.
In particular, Condition (\textasteriskcentered) contains
\eqref{eqdifference}.
\end{remark}

\begin{proof}[Proof of Theorem \ref{th:necessity}]We argue by
induction on $n$. For $n=1$ the assertion is trivial and for $n=2$
we refer to Theorem \ref{thm:twodecomp} and Remark
\ref{rem:symmass}. Let $n>2$ and assume that the statement of the
theorem is true for all $n'<n$. Suppose that for the partition
$B_1\cup B_2\cup\cdots \cup B_N=\{1,2,\dots,n\}$, for $h_j\in
B_j$, $k_j,l_i,l_i'\in \NN$, $(j=1,\dots, N$, $i\in B_j\setminus
\{h_j\}$) and for $z\in A$ the conditions in \eqref{eq:mcond}
hold. We need to prove \eqref{eq:mconclusion}.

First, suppose that there are at least two non-empty blocks in the
partition, and, say, $h_1\in B_1$. We will apply the induction
hypothesis in the following situation. We define $A'$ as the orbit
of $z$ under $B_2\cup\dots\cup B_N$ and
$f':=(\Diff_{T^{k_1}_{h_1}}f)|_{A'}$. Since $f$ has a $(T_1,\dots,
T_n)$-decomposition, the function $f'$ has an invariant
decomposition with respect to the transformations belonging to
$B_2\cup\cdots\cup B_N$. Indeed, by assumption $f=f_1+\cdots +f_n$
with $f_i$ being $T_i$-invariant. We show that
$(\Diff_{T_{h_1}^{k_1}}f_i)|_{A'}=0$ for $i\in B_1$. This is obvious
for $i=h_1\in B_1$, so let us assume $i\in B_1\setminus\{h_1\}$.
Then by \eqref{eq:mcond} for $x\in A'$ we have $x=Sz$ with a
suitable product $S$ of mappings $T_j$, $j\in B_2\cup\dots \cup B_N$
that
\begin{align*}
\Diff_{T_{h_1}^{k_1}}f_i(x)&=f_i(T_{h_1}^{k_1} x)-f_i(x)
=f_i(T_{h_1}^{k_1}Sz)-f_i(Sz)\\
&=f_i(T_{h_1}^{k_1}
T_i^{l_i}Sz)-f_i(T_i^{l_i'}Sz)\\
&=f_i(T_i^{l_i'}Sz)-f_i(T_i^{l_i'}Sz)=0~.
\end{align*}
Since we have the relations
\begin{equation*}
T_{h_j}^{k_j} T_{i}^{l_i}z=T_{i}^{l_i'}z\qquad \mbox{for all $i\in
B_j$, $i\neq h_j$, for all  $j=2,\dots, N$}~,
\end{equation*}
the induction hypothesis gives
\begin{equation*}
\Diff_{T_{h_2}^{k_2}}\dots
\Diff_{T_{h_{N}}^{k_N}}\Diff_{T_{h_1}^{k_1}}f(z)=\Diff_{T_{h_2}^{k_2}}\dots
\Diff_{T_{h_{N}}^{k_N}}f'(z)=0~,
\end{equation*}
hence the assertion.

Second, we suppose that there is only one block in the partition,
i.e., $B_1=\{1,\dots,n\}$. We have to show $\Diff_{T_{h_1}^{k_1}}
f(z)=0$. We can suppose without loss of generality that $h_1=1$ (we
also write $k=k_1$) so \eqref{eq:mcond} becomes
\begin{equation*}
T_1^{k} T_{j}^{l_j}z=T_{j}^{l_j'}z,\qquad \mbox{$j=2,\dots,n$}~.
\end{equation*}
By the condition of decomposability, i.e., the validity of \eqref{invaridecdef}, it suffices to
show only that
\begin{equation*}
f_j(T_1^{k}z)=f_j(z)\qquad \mbox{if}\quad \Diff_{T_j}f_j\eqegy0~.
\end{equation*}
But, indeed, this holds by
\begin{equation*}
f_j(T_1^{k}z)=f_j(T_1^{k}T_j^{l_j}z)=f_j(T_j^{l'_j}z)=f_j(z)~.
\end{equation*}
\end{proof}

 We thank to Tam\'{a}s Keleti for suggesting the above
formulation of the Condition (\textasteriskcentered). This allows
for the following nice reformulation in Abelian groups, which was
also suggested by him.

\begin{remark}\label{groupcondi} It is particulary interesting to
formulate Condition (\textasteriskcentered) in the special case,
when $A$ is an Abelian group and the transformations $T_i$ are
translations by $a_i\in A$, i.e., $T_i x=x+a_i$ for $x\in A$,
$i=1,\dots,n$. In this case, Condition (\textasteriskcentered) takes
the following form.

\vskip1em\noindent\hbox
to\hsize{\vtop{\hsize=0.1\hsize\vfill\noindent{\rm
(\textasteriskcentered\textasteriskcentered)}\vfill}\vtop{\hsize=0.9\hsize
\noindent\begin{minipage}{\hsize} Suppose that whenever for a
partition $B_1\cup B_2\cup\cdots \cup B_N=\{1,2,\dots,n\}$ with
distinguished elements $h_j\in B_j$ and for the natural numbers
$k_j\in \NN$, $j=1,\dots, N$
\begin{equation}\label{eq:mcond2}
a_i\mbox{ divides$^1$ }k_j\cdot a_{h_j}\qquad \mbox{for all $i\in
B_j\setminus\{h_j\}$ and for all $j=1,\dots, N$}
\end{equation}
then,
\begin{equation}\label{eq:mconclusion2}
\Diff_{T_{h_1}^{k_1}}\dots \Diff_{T_{h_{N}}^{k_N}}f=0
\end{equation}
must also be satisfied.
\end{minipage}}}\vskip1em
\footnotetext[1]{For $a,b\in A$, we say that $a$ \emph{divides} $b$
if there is $n\in \NN$ with $n\cdot a=b$.}
\end{remark}

We saw in Theorem \ref{thm:twodecomp} that Condition (\textasteriskcentered) -- or even a
subset of the conditions listed in it -- provides also a sufficient condition if $n=2$. Next we
push this further.

\begin{theorem}\label{thm:three}
Suppose that $T_1=T$, $T_2=S$ and $T_3=U$ commute and the function $f$ satisfies Condition
(\textasteriskcentered). Then $f$ has a $(T,S,U)$-invariant decomposition.
\end{theorem}
The proof will be based on the following series of lemmas.

\begin{lemma}\label{lem:solvg}
Let $g\in\RR^A$ be a function and $T$ be a transformation of the set $A$.  The following
statements are equivalent.
\begin{enumerate}[i)]
\item There exists a function $h$ for which $\Diff_T h=g$.
\item $\sum_{i=0}^{k-1}g(T^i x)=0$ whenever $T^kx=x$, $x\in A$, $k\in\NN$.
\end{enumerate}
\end{lemma}
\begin{proof}
If i) holds and $T^kx=x$ is satisfied for some $x\in A$ and $k\in\NN$, then
$$
\sum_{i=0}^{k-1}g(T^i x)=\sum_{i=0}^{k-1}\Bigl(h(T^{i+1} x)-h(T^ix)\Bigr)=h(T^{k}x)-h(x)=0.
$$

Suppose now that ii) holds. We define the equivalence relation: $x\ekv y$ iff $T^k x=T^l y$
with some $k,l \in\NN$. By the axiom of choice, we select a representative of each equivalence
class. Then it suffices to give a proper construction of $h$ on an arbitrarily given
equivalence class $B$ with representative $x_0$, say.

If both
\begin{equation}\label{twoekviv}
T^m x = T^n x_0 \qquad \hbox{and} \qquad T^{m'} x = T^{n'} x_0~,
\end{equation}
then also
\begin{equation}\label{allekvforms}
T^{m+n'} x = T^{n+n'} x_0 = T^{n+m'} x~.
\end{equation}
First suppose that $m'\geq m$ and $n'\geq n$, and compute with $z:=T^m x = T^n x_0$ and
$M:=\min(n'-n,m'-m)$ and $N:=\max(n'-n,m'-m)$ the relations
\begin{align}\label{geqations}
& \sum_{i=0}^{n'-1}g(T^i x_0) -  \sum_{j=0}^{m'-1} g(T^j x) -
\sum_{i=0}^{n-1} g(T^i x_0) + \sum_{j=0}^{m-1}g(T^j x)
\notag \\
&\qquad \qquad =\sum_{i=n}^{n'-1} g(T^i x_0) - \sum_{j=m}^{m'-1} g(T^j x)
\notag\\
&\qquad \qquad=\sum_{i=0}^{n'-n-1} g(T^i z) - \sum_{j=0}^{m'-m-1} g(T^j z) = \pm
\sum_{l=M}^{N-1} g(T^l z)~.
\end{align}
Now suppose, e.g., that $N=n'-n\geq M=m'-m$ (the opposite case being similar). Because of
\eqref{twoekviv} we have
$$
T^M z=T^{m'-m} z = T^{m'} x = T^{n'}x_0=T^{n'-n}z=T^N z ~,
$$
which, in view of ii), immediately gives the vanishing of \eqref{geqations}.

In case we do not have both the conditions $m'>m$ and $n'>n$ let us take $m'':=m+m'+n'$,
$n'':=n+n'+m'$, and, taking also \eqref{allekvforms} into account, apply the known case to
$n,m$ and $n'',m''$ as well as to $n',m'$ and $n'',m''$ separately. These considerations then
tell us that \eqref{geqations} is always $0$ and so
\begin{equation}\label{hdefeq}
h(x):=\sum_{i=0}^{n-1} g(T^i x_0) - \sum_{j=0}^{m-1}g(T^j x)~, \quad \hbox{whenever} \quad T^m
x=T^n x_0
\end{equation}
is a correct definition of a function $h$. In $T^m x= T^n x_0$ we
may suppose that $m>0$, so by $T^{m-1}(Tx)=T^n x_0$ and using
\eqref{hdefeq} we have
\begin{equation*}
h(Tx)-h(x)=\sum_{i=0}^{n-1} g(T^i x_0) - \sum_{j=0}^{m-2}g(T^j T
x)-\sum_{i=0}^{n-1} g(T^i x_0) + \sum_{j=0}^{m-1}g(T^j x)=g(x),
\end{equation*}
hence the assertion follows.
\end{proof}

\begin{lemma}\label{lem:persolv}If $G:A\to \RR$ is an arbitrary function and $T:A\to A$ is an arbitrary
transformation, then there is a function $g$ and a $T$-invariant function $\gamma$ such that
$$
\Diff_T g=G+\gamma~.
$$
\end{lemma}
\begin{proof}
First of all we define an equivalence relation $\ekv$. We say that $x$ and $y$ are equivalent,
$x\ekv y$, if there exists $n,m\in\NN$ such that $T^ny=T^mx$. Of course this is indeed an
equivalence relation, and the equivalence class of an element $x\in A$ is denoted by $B_x$. In
view of Lemma \ref{lem:solvg} we are looking for some $T$-invariant $\gamma$ satisfying
\begin{equation*}%\label{eq:eq3}
\sum_{i=0}^{k-1}\gamma(T^i x)=-\sum_{i=0}^{k-1}G(T^i x)~,\qquad\mbox{whenever} \quad T^kx=x,
~k\in\NN~.
\end{equation*}
By $\Diff_T\gamma\eqegy0$ this is equivalent to the assertion that for every equivalence class
$B_z$ there is a constant $\gamma=\gamma(B_z)$ such that
$$
-\gamma=\frac{1}{k}\sum_{i=0}^{k-1}G(T^i x)~, \qquad \mbox{if} \quad T^k x=x,~ k\in\NN^+, ~x\in
B_z~.
$$
Suppose that $y\ekv x$, $T^k x=x$ and $T^l y=y$, $k,l\in\NN^+$. By $x,y\in B_z$ we have
$T^ax=T^b y$ for some $a,b\in\NN$, and for $K=kl$ the equations $T^Kx=x$, $T^K y=y$ hold. Now
\begin{equation}\label{eq:kl}
\begin{split}
\frac{1}{l}\sum_{i=0}^{l-1}G(T^i y)&=\frac{1}{K}\sum_{i=0}^{K-1}G(T^i
y)=\frac{1}{K}\sum_{i=b}^{K-1+b}G(T^i
y)\\
&=\frac{1}{K}\sum_{i=0}^{K-1}G(T^iT^ax)=\frac{1}{k}\sum_{i=0}^{k-1}G(T^i x)~,
\end{split}
\end{equation}
which means that this quantity is constant for $x,y$ with the
above properties. Define
$$
\gamma(x):=\frac{1}{k-l}\sum_{i=l}^{k-1}G(T^i x)~, \qquad \mbox{if} \quad T^k x=T^l x,~
k,l\in\NN, \:k>l~,
$$
and $\gamma(x)$ arbitrary if there are no such $k,l$. By the above
argument this definition is independent of the particular choice of
$k,l$. If $x\ekv y$ and for $x$ there are no $k,l$ satisfying $T^k
x=T^l x$, neither can exist such $k,l$ for $y$. Thus we see that
$\gamma$ can be chosen to be constant on $B_z$. The proof is hence
complete.
\end{proof}

\begin{lemma}\label{lem:twopersolv} Let $T$ and $S$ be commuting transformations of the set $A$, and let $G:A\to \RR$ be an
$S$-invariant function. Then there exist functions $\gamma$ and $g$ such that $\Diff_S
\gamma=\Diff_S g\eqegy0$, $\Diff_T \gamma\eqegy0$ and
$$
\Diff_T g=G+\gamma~.
$$
\end{lemma}
\begin{proof}
Again we define an equivalence relation, $x\ekv y$ if $S^nx=S^my$ for some $n,m\in\NN$. Because
of commutativity $Tx\ekv Ty$ whenever $x\ekv y$. Let us consider
$\widetilde{A}:=A/\hskip-5pt\ekv$ and define $\widetilde{T}(B_x):=B_{Tx}$, where, in general,
$B_z$ stands for the equivalence class of $z$. By the above the transformation $\widetilde{T}$
is well-defined. Since $G$ is $S$-invariant, it is constant on each equivalence class $B_x$, so
$\widetilde{G}(B_x):=G(x)$ is a correct definition of a real-valued function on $\widetilde{A}$.
Applying Lemma \ref{lem:persolv} to $\widetilde{A}$, $\widetilde{T}$, $\widetilde{G}$ we obtain
the functions $\widetilde{\gamma}$ and $\widetilde{g}$ with $\Diff_{\widetilde{T}}
\widetilde{g}=\widetilde{G}+\widetilde{\gamma}$ and $\Diff_{\widetilde{T}} \widetilde{\gamma}
\eqegy 0$. Defining $g$ and $\gamma$ to be constant on each equivalence class of $\ekv$:
$$
g(x):=\widetilde{g}(B_x),\qquad \gamma(x)=\widetilde{\gamma}(B_x)
$$
we see immediately that $\Diff_S \gamma=\Diff_S g\eqegy 0$. Obviously $\Diff_{\widetilde{T}}
\widetilde{g}=\widetilde{G}+\widetilde{\gamma}$ implies $\Diff_T g=G+\gamma$ and
$\Diff_{\widetilde{T}} \widetilde{\gamma} \eqegy 0$ implies $\Diff_T \gamma \eqegy 0$. This
completes the proof.
\end{proof}

\begin{lemma}\label{lem:diffinvsolv} Let $T,S$ be commuting transformations of $A$ and let
$G:A\to\RR$ be a function satisfying $\Diff_S G=0$. Then there
exists a function $g: A\to \RR$ satisfying both $\Diff_S g= 0$ and
$\Diff_T g =G$ if and only if
\begin{equation}\label{twotransfcondi}
\sum_{i=0}^{k-1} G(T^i x) =0 \qquad\hbox{whenever}\quad T^kS^lx=S^{l'}x,\,x\in A,\, k,l,l'\in
\NN~.
\end{equation}
\end{lemma}
\begin{proof} It is easy to check that existence of a function
$g$ with the above requirements implies \eqref{twotransfcondi} (cf Lemma \ref{lem:solvg}),
hence we are to prove sufficiency of this condition only.

We can argue similarly as in the proof of Lemma \ref{lem:twopersolv}, but using Lemma
\ref{lem:solvg} in place of Lemma \ref{lem:persolv}.

Namely, as $G$ is $S$-invariant, we can consider the equivalence relation $x\ekv y$ iff $S^n
x\ekv S^m y$ for some $n,m\in\NN$, and define $\widetilde{G}$ on the set of equivalence classes
$\widetilde{A}:=A/\hskip-5pt\ekv$ as the common value $G(x)$ of the function $G$ on the whole
class $B_x$. Also, by commutativity, $T$ generates a well-defined transformation
$\widetilde{T}$ of $\widetilde{A}$. With this definition, Lemma \ref{lem:solvg} applies to
$\widetilde{G}$ and $\widetilde{T}$; note that in $\widetilde{A}$ two classes are related with
respect to $\widetilde{T}$ as $\widetilde{T}^{k'}(B_x)=B_x$ iff the condition in
\eqref{twotransfcondi} holds. Therefore, \eqref{twotransfcondi} is equivalent to condition ii)
of Lemma \ref{lem:solvg} when applied to the function $\widetilde{G}$ on the set
$\widetilde{A}$ and the transformation $\widetilde{T}$. The ``lift up'' $g$ of $\Tilde{g}$ as
in the proof of  Lemma \ref{lem:twopersolv} will be appropriate.
\end{proof}

\begin{remark}\label{r:gammaconst} Combining the last two
lemmas, one can see that $\Diff_T g=G$ is equivalent to the requirement that $\gamma(x)=0$
whenever $ T^k S^lx=S^{l'}x,\,x\in A$ and $k,l,l'\in \NN$; moreover, any proper $\gamma$ in
Lemma \ref{lem:twopersolv} must satisfy
\begin{equation}\label{gammaformula}
\gamma(x)=-\frac 1k \sum_{i=0}^{k-1} G(T^i x) \qquad\hbox{whenever}\quad T^kS^lx=S^{l'}x,\,x\in
A,\, k\in\NN^+,l,l'\in \NN\,.
\end{equation}
In particular,
\begin{equation}\label{gammaformula2}
\gamma(x)=-\frac 1{k-k'} \sum_{i=k'}^{k-1} G(T^i x) \quad\hbox{if}\quad T^kS^lx=T^{k'}S^{l'}x,
k,k',l,l'\in \NN,k-k'>0~.
\end{equation}
Furthermore, looking at the proof of Lemma \ref{lem:persolv}, we see that if no such conditions
as in \eqref{gammaformula2} are satisfied for $x$, then $\gamma(x)$ can be chosen to be an
arbitrary constant on each equivalence class of $\ekv_{TS}$, where $x\ekv_{TS} y$ if $S^a T^b
x= S^{a'} T^{b'} y$, for some $a,b,a',b'\in\NN$.
%\todo{Ez itt nem\\ vilagos.\\ Hogy is\\ vesszuk $\gamma$-t?\\
%can not be arbitrary}On the other hand, it is easy to see that once we define a function
%$\gamma$ according to this requirement, take $\gamma(y):=\gamma(x)$ for any $y$ with $S^a T^b
%x= S^{a'} T^{b'} y$ and arbitrarily otherwise, then we get a solution of the requirements of
%Lemma \ref{lem:twopersolv}.
\end{remark}

\begin{proof}[Proof of Theorem \ref{thm:three}]
Take $\Diff_T f=F$. By taking $B_1=\{1\}$, $B_2=\{2\}$ and
$B_3=\{3\}$  in Condition (\textasteriskcentered) we have
\begin{equation}\label{eq:threeper}
\Diff_T\Diff_S\Diff_U f\eqegy0,\qquad\mbox{that is} \quad
\Diff_S\Diff_U F\eqegy0~.
\end{equation}
Further, also by Condition (\textasteriskcentered) for the partition
$B_1=\{1\}$, $B_2=\{2,3\}$, $h_2=3$, if
$U^{k+k'}S^nx=U^{k'}S^{n'}x$, then
\begin{equation}\label{eq:orbs}
\Diff_T\Diff_{U^k}f(U^{k'}x)=0,\qquad\mbox{that is} \quad F(U^{k+k'}x)=F(U^{k'}x)~.
\end{equation}
The equations \eqref{eq:threeper} and \eqref{eq:orbs} show that \ref{enum:1:ii}) of Theorem
\ref{thm:twodecomp} is fulfilled. This implies the existence of $S$- and $U$-invariant
functions $H$ and $L$ respectively with
$$
F=H+L~.
$$
We apply Lemma \ref{lem:twopersolv} to obtain the real-valued functions $h, l,\chi,\lambda$ with
\begin{align}\label{alldiffeq}
\Diff_T h& =H+\chi, &\Diff_Sh&=\Diff_S\chi=\Diff_T\chi\eqegy0~,
\\
\Diff_T l& =L+\lambda,& \Diff_Ul&=\Diff_U\lambda=\Diff_T\lambda\eqegy0~.\notag
\end{align}
Define $g:=f-h-l$, then $f=g+h+l$ and $\Diff_S h = \Diff_U l = 0$, while $\Diff_T g= \Diff_T
(f-h-l) = F - \Diff_T h - \Diff_T l$. So using now the decomposition $F=H+L$ and
\eqref{alldiffeq}, we arrive at $\Diff_T g= (H+L) - (H+\chi) - (L+\lambda)= - \chi- \lambda$.

To illustrate the merit of the next argument, let us assume first
that $T^kx=x$ for some $x\in A$ and $k\in \NN^+$. Then we can refer
to Lemma \ref{lem:solvg}. We have seen that the function
$\gamma:=-(\chi+\lambda)=\Diff_T g$, hence condition i) of the lemma
is satisfied and $\gamma$ must satisfy the equivalent condition ii).
On the other hand, $\gamma$ is also $T$-invariant by construction
(see \eqref{alldiffeq}), hence ii) of Lemma \ref{lem:solvg} can be
satisfied if only $\gamma(T^i x)=0$ for all $i=0,\dots,k$.
Therefore, in case $T^kx=x$ for some $x\in A$ and $k\in \NN^+$, we
already have $\Diff_T g (x)=\gamma(x) =0$. Note also that, because
of the $T$-invariance of $\gamma$, for any $y\in A$ with $T^n x =
T^m y$ for some $n,m\in \NN$ one must have $\Diff_T g (y) =\Diff_T g
(x)$, in particular $\Diff_T g (y)=0$ if $x$ is as before.

Our aim is to obtain the same thing in general, for all over $A$. In the definition of $\chi$
and $\lambda$ we may have certain flexibility. Exploiting this and choosing both functions
carefully we will have $\gamma=-(\chi+\lambda)=0$. For this purpose, we define an equivalence
relation
\begin{equation}\label{classdef}
x\ekv y\quad\quad \mbox{iff}\quad\quad T^a S^b U^c x = T^{a'} S^{b'} U^{c'} y\quad\mbox{for
some $a,b,c,a',b',c'\in\NN$}~.
\end{equation}
It suffices to restrict considerations to one equivalence class $B_z=\{y\in A:\:z\ekv y\}$,
so without loss of generality we can work only on $B_z$ assuming tacitly $A=B_z$.

By Remark \ref{r:gammaconst}, we \emph{can not} choose $\chi(x)$ for
some $x$ arbitrarily if some relation of the type
\begin{equation}\label{chicondi}
 T^kS^lx=T^{k'}S^{l'}x,\qquad\quad\,k,k',l,l'\in\NN,\:k>k'
\end{equation}
holds. Let us call such points $x$ $(S,T)$-prescribed.

Suppose first that there are neither  $(S,T)$-prescribed nor
$(U,T)$-prescribed points. In this case recalling Lemma
\ref{lem:diffinvsolv} we can choose both $\chi$ and $\lambda$ to
be, e.g., constant $0$ on $A$.

Next, by symmetry, we can assume that there are e.g.
$(S,T)$-prescribed points. We will show that in this case $\chi$
can be also chosen to be a constant. So let now $x\in A$ be fixed
and satisfying \eqref{chicondi}. Then
\begin{equation}\label{chicondi2}
\chi(x)=-\frac 1{k-k'} \sum_{i=k'}^{k-1} H(T^i x)
\end{equation}
must hold by Remark \ref{r:gammaconst}. Moreover, relations as in
\eqref{chicondi} hold for all $y\in A$, $y\ekv_{TS} x$ (where
$x\ekv_{TS} y$ iff $S^a T^b x= S^{a'} T^{b'} y$ for some
$a,b,a',b'\in\NN$, as in Remark \ref{r:gammaconst}). Conversely, if
$y$ is not $(S,T)$-prescribed, then $y$ can not be in $\ekv_{TS}$
relation with the above $x$, and by Remark \ref{r:gammaconst}, the
value of $\chi$ can be chosen to be arbitrary constant on the whole
$\ekv_{TS}$-class of $y$. So let this constant be $\chi(x)$ ($x$ is
the above fixed element). Now, as there might exist elements $y\in
A$ having \eqref{chicondi}, we show that $\chi(y)=\chi(x)$ for all
such $y$, too, so for all $y\in A$ regardless whether $y\ekv_{TS}x$
holds or not. To this end, notice first that the relation
$T^kS^lUx=T^{k'}S^{l'}Ux$ is also valid. So using $\Diff_U L=0$ and
thus $\Diff_U H= \Diff_U F = \Diff_T \Diff_U f$, we compute
\begin{align}
\Diff_U \chi(x) & = \chi(Ux)-\chi(x)=-\frac 1{k-k'} \sum_{i=k'}^{k-1}
H(T^i Ux) + \frac 1{k-k'} \sum_{i=k'}^{k-1} H(T^i x) \notag \\
&= -\frac 1{k-k'} \sum_{i=k'}^{k-1} \Diff_U H(T^i x) = -\frac 1{k-k'} \sum_{i=k'}^{k-1} \Diff_U
F (T^i x)\notag\\
&=\Diff_U \Diff_{T^{k-k'}} f(T^{k'}x) =0~,\label{diffchi}
\end{align}
 the last step being an application of Condition
(\textasteriskcentered) for the partition $B_1=\{3\}$,
$B_2=\{1,2\}$, $h_2=1$ and the equation $T^kS^lx=T^{k'}S^{l'}x$
(remember the same argument works for $y$, too). Because by
assumption $T^a S^b U^c x = T^{a'} S^{b'} U^{c'} y$ holds for some
$a,b,c,a',b',c'\in\NN$, we obtain (using also $\chi(Ux)=\chi(x)$,
$\chi(Uy)=\chi(y)$, as proved above, and the $S$- and $T$-invariance
of $\chi$) that
$$
\chi (x)=\chi(T^a S^b U^c x)=\chi(T^{a'} S^{b'} U^{c'} y)=\chi(y),
$$
which shows that $\chi$ is indeed constant on the whole of $A$.

Now, if there is no $(U,T)$-prescribed point, then there is absolute
freedom in choosing $\lambda$ (as long as it is constant on the
equivalence classes), so we can define it to be the negative of the
constant value of $\chi$. We obtain $\gamma=-(\chi+\lambda)=0$ as
required.

Finally, suppose that there are both $(S,T)$-prescribed and
$(U,T)$-prescribed points, say $x$ and $x'$, which then satisfy
\eqref{chicondi} and
\begin{equation}\label{lcondi}
T^mU^n x'=T^{m'}U^{n'}x'\qquad\mbox{for
some}\quad\,m,m',n,n'\in\NN^+,m>m'~,
\end{equation}
respectively. Since  $x\ekv x'$, by definition $T^a S^b U^c x =
T^{a'} S^{b'} U^{c'} x'(=:y)$ holds for some $a,b,c,a',b',c'\in\NN$.
But then $y$ is both $(S,T)$- and $(U,T)$-prescribed, that is both
\eqref{chicondi} and \eqref{lcondi} holds with $y$ replacing $x$ and
$x'$, respectively. From this we conclude that
\begin{equation}\label{eq:combi}
\begin{split}
T^{(m-m')(k-k')}U^{n(k-k')}T^{m'+k'}y&=U^{n'(k-k')}T^{m'+k'}y~,\quad\mbox{and}\\
T^{(m-m')(k-k')}S^{l(m-m')}T^{m'+k'}y&=S^{l'(m-m')}T^{m'+k'}y~.
\end{split}
\end{equation}
To see the first identity here, we can write, by \eqref{lcondi}
(with $x'$ replaced by $y$) and applying $T^{k'}$ to both sides,
that
\begin{equation*}
 T^{m-m'}U^{n}T^{m'+k'}y=U^{n'}T^{m'+k'}y~.
\end{equation*}
This is just the required equation provided $k-k'=1$, and the
general case follows from this inductively. The second line in
\eqref{eq:combi} is proved analogously.

Recall that in case of presence of $(S,T)$-prescribed points the
function $\chi$ can be chosen having a constant value. The same
applies for $\lambda$ in case there exists some $(U,T)$-prescribed
point. So let us now assume that both functions are defined as
constants all over $A$. Then it remains to show that these
constant functions sum to 0 at some, hence on all points of $A$.

By \eqref{eq:combi} and writing $z:=T^{m'+k'}y$, we arrive at
$T^{K}U^Nz =U^{N'}z$ and $T^{K}S^Lz=S^{L'}z$, with appropriate
$K,L,N,L',N'\in\NN$, $K>0$. Thus, by Remark \ref{r:gammaconst}, we
must have
\begin{align}\label{chilambda}
  \chi(z)&= - \frac 1{K} \sum_{i=0}^{K-1} H(T^i z) ~,\quad\mbox{since}\quad T^{K}S^Lz=S^{L'}z~,\quad\mbox{and}\\
  \lambda(z)&= - \frac 1{K} \sum_{i=0}^{K-1} L(T^i z)~,\quad\mbox{since}\quad T^{K}U^Nz=U^{N'}z~.\notag
\end{align}
Summing these and using the decomposition of $F$ we obtain
\begin{equation}\label{gammanull}
\gamma(z)=- (\chi(z)+\lambda(z))= \frac 1{K} \sum_{i=0}^{K-1}
F(T^i z) = \frac 1{K} \Diff_{T^{K}} f (z) =0~,
\end{equation}
by Condition (\textasteriskcentered) for the partition
$B_1=\{1,2,3\}$ with $h_1=1$ and in view of the equations on the
right of \eqref{chilambda}. That is, $\chi+\lambda$ is zero, and
the proof is complete.
\end{proof}

\begin{question}
 We close this section by the natural question if Condition (\textasteriskcentered) is
equivalent to \eqref{invaridecdef} for all $n\in\NN$ ($n>3$).
\end{question}

\section{Further results for unrelated transformations}

We call two commuting transformations $S,T$ on $A$
\emph{unrelated}, iff $T^nS^k x= T^m S^l x$ can occur only if
$n=m$ and $k=l$. In particular, then neither of the two
transformations can have any cycles in their orbits, nor do their
joint orbits have any recurrence.

If all pairs among the transformations $T_j$ ($j=1,\dots,n$) are
unrelated, then Condition (\textasteriskcentered) degenerates, as in
\eqref{eq:mcond} we necessarily have that all blocks $B_j$ are
singletons. We saw in Remark \ref{mzerocondi} that it is exactly the
difference equation \eqref{eqdifference}.

As an application in a special situation, consider now the case
when the set $A:=\RR$ and the transformations are just shifts by
real numbers. It is easy to see that $T_{\alpha}$ and
$T_{\beta}$, the shift operators by $\alpha\in\RR$ and $\beta
\in\RR$, are unrelated iff $\alpha/\beta$ is irrational.
Therefore, for $n=3$ we obtain the following special case from
Theorem \ref{thm:three}.

\begin{corollary}\label{cr:irrational} Let $\alpha_i$ ($i=1,\dots,n$)
be nonzero real numbers so that $\alpha_i/\alpha_j$ are irrational whenever $1\le i\ne j \le
n$. Then the conditions \eqref{eq:periodd} and \eqref{eq:kern} are equivalent.
\end{corollary}

We stated the above corollary for general $n$ since for unrelated
transformations it can be proved for any $n \in \NN$. In fact,
the following more general form holds.

\begin{theorem}\label{thm:nomix}
If the transformations $T_j$ ($j=1,\dots,n$) are pairwise
(commuting and) unrelated, then the difference equation
\eqref{eqdifference} is equivalent to the existence of some
invariant decomposition \eqref{invaridecdef}.
\end{theorem}

\begin{proof} We argue by induction. The cases of small $n$ are obvious.
Existence of an invariant decomposition \eqref{invaridecdef} clearly implies the difference
equation \eqref{eqdifference} for any set of pairwise commuting transformations, unrelated or
not, hence it suffices to deal with the converse direction.

Let $F:=\Diff_{T_{n+1}} f$. As the $n+1$-level difference equation of $f$ is inherited by $F$
as an $n$-level one, by the inductive hypothesis we can find an invariant decomposition of $F$
in the form
\begin{equation}\label{Fdecomp}
F=F_1+\cdots +F_n \,,\qquad\text{where}\qquad  \Diff_{T_j} F_j \eqegy0 \quad (j=1,\dots,n)\,.
\end{equation}
Since $T_{n+1}$ and $T_j$ are unrelated for $j=1,\dots,n$, the
condition \eqref{twotransfcondi} in Lemma \ref{lem:diffinvsolv}
is void, and therefore the "lift ups" $f_j$ with $\Diff_{T_j}
f_j=0$, $\Diff_{T_{n+1}}f_j=F_j$ exist for all $j=1,\dots,n$.
Therefore, $f_{n+1}:=f-f_1-\dots-f_n$ provides a function
satisfying $\Diff_{T_{n+1}} f_{n+1} = F - F_1-\dots - F_n =0$,
whence a decomposition of $f$ is established.
\end{proof}

\section{On invariant decompositions of bounded functions}

Finally, let us mention a complementary result, which concerns
bounded functions, thus is not fully in scope here, but is similar
in nature regarding the absolutely unrestricted structural
framework of transformations and functions.

\begin{proposition} Let $A$ be any set, $T,S:A\to A$ arbitrary
commuting transformations, and let $G:A\to \RR$ be any function satisfying $\Diff_S G=0$. Then
the following two assertions are equivalent.
\begin{enumerate}[i)]
\item$\exists\:\: H: A\to \RR$ bounded function such that $\Diff_T H=G$
and $\Diff_S H= 0$. \\
\item $\exists\:\: C<+\infty$ constant such that $\left|\sum_{i=1}^{m-1}
G(T^i x) \right| \leq C $ whenever $x\in A$ and $m\in \NN$.
\end{enumerate}
Moreover, one has the relations $C\leq \|H\|_{\infty}\leq 2C$.
\end{proposition}
\begin{proof}
The implication i)$ \Rightarrow$ ii) is immediate with $C:=2\|H\|_{\infty}$, since
$$
\sum_{i=1}^{m-1} G(T^i x)  = \sum_{i=1}^{m-1} \Diff_T H(T^i x) = H(T^m x)-H(x)\,.
$$

The proof of the converse direction ii) $\Rightarrow$ i) goes along similar lines to the
above, hence we skip the details.
\end{proof}

We mention this as an example of the case when the class of functions on $A$ we deal with is
$B(A)$, the set of all bounded functions. It is known that $B(A)$ has the decomposition
property, see \cite{gajda:1992} and \cite{laczkovich/revesz:1990}, but the exact norm
inequalities are not known and very likely depend on the transformations, in particular
properties like unrelated and alike. On the other hand it is remarkable, that if $f$ is
bounded, then no further conditions, neither on the transformations nor on $f$ are involved:
\eqref{eqdifference} itself implies \eqref{invaridecdef}. It would be interesting, but perhaps
difficult, to determine the best general bound for the norms of individual terms in
\eqref{invaridecdef} once $\|f\|$ is given.

\section*{Acknowledgment}

We thank  Tam\'{a}s Keleti as well as the anonymous referee for their
constructive comments and suggestions

\noindent B.~Farkas\\
Technische Universit\"{a}t Darmstadt\\
Fachbereich Mathematik, AG4\\
Schlo\ss{}gartenstra\ss{}e 7, D-64289, Darmstadt, Germany\\
\texttt{farkas@mathematik.tu-darmstadt.de}

 \vskip2ex
\noindent Sz.~Gy.~R\'{e}v\'{e}sz\\
Alfr\'{e}d R\'{e}nyi Institute of Mathematics\\
Hungarian Academy of Sciences\\ Re\'{a}ltanoda utca 13--15, H-1053, Budapest, Hungary\\
\texttt{revesz@renyi.hu}\\[1ex]
and\\[1ex] Institut Henri Poincar\'{e} \\
11 rue Pierre et Marie Curie,\\
75005 Paris, France \\
\texttt{Szilard.Revesz@ihp.jussieu.fr}
\end{document}